\documentclass{amsart}
\usepackage{graphicx}
\usepackage{amscd}

\newcommand{\comment}[1]{\marginpar{\sffamily{\noindent\tiny #1
   \par}\normalfont}}
\renewcommand{\comment}[1]{}

\newcommand{\C}{{\mathbb C}}
\newcommand{\Z}{{\mathbb Z}}
\newcommand{\R}{{\mathbb R}}

\newcommand{\co}{\colon}

\newtheorem{theorem}{Theorem}[section]
\newtheorem{lemma}[theorem]{Lemma}
\newtheorem{proposition}[theorem]{Proposition}

\theoremstyle{definition}
\newtheorem{definition}[theorem]{Definition}
\newtheorem{example}[theorem]{Example}
\newtheorem{remark}[theorem]{Remark}

\begin{document}
\title[$\mu$-constancy does not imply constant bi-Lipschitz type]
{$\mu$-constancy does not imply constant\\ bi-Lipschitz type}
\author{Lev Birbrair}
\address{Departamento de Matem\'atica, Universidade Federal do Cear\'a
(UFC), Campus do Pici, Bloco 914, Cep. 60455-760. Fortaleza-Ce,
Brasil} \email{birb@ufc.br}
\author{Alexandre Fernandes}
\address{Departamento de Matem\'atica, Universidade Federal do Cear\'a
(UFC), Campus do Pici, Bloco 914, Cep. 60455-760. Fortaleza-Ce,
Brasil} \email{alexandre.fernandes@ufc.br}
\author{Walter D.
  Neumann} 
\address{Department of Mathematics, Barnard College,
  Columbia University, New York, NY 10027}
\email{neumann@math.columbia.edu}


\begin{abstract}
  We show that a family of isolated complex hypersurface singularities
  with constant Milnor number may fail, in the strongest sense, to
  have constant bi-Lipschitz type. Our example is the 
  Brian\c con--Speder family $X_t:=\{(x,y,z)\in\C^3 ~|~
  x^5+z^{15}+y^7z+txy^6=0 \}$ of normal complex surface germs; we show
  the germ $(X_0, 0)$ is not bi-Lipschitz homeomorphic with respect to
  the inner metric to the germ $(X_t,0)$ for $t\ne 0$.
\end{abstract}

\maketitle

\section{Introduction}

Given a germ $(X,p)$ of a point of a complex analytic set, a choice of
generators $x_1,\dots,x_N$ of its local ring gives an embedding of
$(X,p)$ into $(\C^N,0)$. It then carries two induced metric space
structures: the ``outer metric'' induced from distance in $\C^N$ and
the ``inner metric'' induced by arc-length of curves on $X$.  In the
Lipschitz category each of these metrics is independent of choice of
embedding: different choices give metrics for which the identity
map is a bi-Lipschitz homeomorphism.  These metric structures have so
far seen much more study in real algebraic geometry than in the
complex algebraic world. The inner metric, which is given by a
Riemannian metric off the singular set, is the one that interests us
most here. It is determined by the outer metric, so germs that are
distinguished by their inner metrics are certainly distinguished by
their outer ones.

It is easy to see that two complex curve germs with the same number of
components are bi-Lipschitz equivalent (inner metric). So for curve
germs bi-Lipschitz geometry is equivalent to topology. This is even so
for outer bi-Lipschitz geometry of plane curves: two germs of
algebraic curves in $\C^2$ are bi-Lipschitz homeomorphic for the outer
metric if and only if they are topologically equivalent as embedded
germs (Teissier--Pham \cite{TP}, Fernandes \cite{F})

We show here that for complex surface germs the picture is very
different. Our main result (announced in \cite{BFN2}) is that in the
Brian\c con-Speder family \cite{BS}
$$X_t:=\{(x,y,z)\in\C^3 ~|~ x^5+z^{15}+y^7z+txy^6=0 \}$$
of germs of algebraic complex surfaces in $\C^3$, which are of
constant embedded topological type, $X_0$ and $X_t$ for $t\ne 0$ are
not bi-Lipschitz homeomorphic, even for the inner metric. In particular,
$\mu$-constancy does not imply bi-Lipschitz equisingularity.

The criterion we use is the existence of so-called \emph{separating
  sets}. A separating set of a complex surface germ $(X,p)$ is a
``thin" 3-dimensional subset $Y$ through $p$ which separates $X$ into
two ``fat" subsets (precise definitions are below). The existence or
non-existence of separating sets is a bi-Lipschitz invariant, and it
turns out that the existence of a separating set is a fairly common
thing in the geometry of complex algebraic surfaces. For example,
$A_k$-singularities were shown to have separating sets for $k>1$ and odd
in \cite{BF} (Theorem \ref{th:sep} below is more
general). We prove that in the Brian\c con-Speder family $X_t$ has a
separating set when $t\neq 0$ and does not have a separating set when
$t=0$.
\subsection*{Acknowledgements}
The authors acknowledge research support under the grants: CNPq grant
no 300985/93-2 (Lev Birbriar),  CNPq
grant no 300393/2005-9 (Alexandre Fernandes) and  NSF grant no.\
  DMS-0456227 (Walter Neumann). They are greatful for useful
  correspondence with Bruce Kleiner and Frank Morgan. 
\section{Separating sets}

\begin{figure}[ht]
    \centering
\includegraphics[width=.33\hsize]{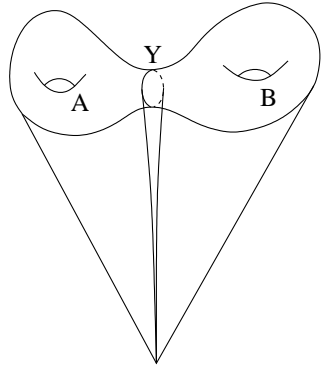}
    \caption{Separating set}
    \label{Figure 1}
  \end{figure}

Let $V\subset\R^n$ be a $k$-dimensional rectifiable subset. Recall
that the inferior and superior $k$--densities of $V$ at the point
$x_0\in\R^n$ are defined by:
$$\underline\Theta^k(V,x_0) = \lim_{\epsilon\to 0^+}
\inf\frac{\mathcal{H}^k(V\cap \epsilon B(x_0))}{\eta\epsilon^k}$$
and
$$\overline\Theta^k(V,x_0)=\lim_{\epsilon\to 0^+}
\sup\frac{\mathcal{H}^k(V\cap \epsilon B(x_0))}{\eta\epsilon^k}\,,$$
where $\epsilon B(x_0)$ is the $n$--dimensional ball of radius
$\epsilon$ centered at $x_0$, $\eta$ is the volume of the
$k$--dimensional unit ball and $\mathcal H^k$ is
$k$--dimensional Hausdorff measure in $\R^n$. If
$$\underline\Theta^k(V,x_0)=\theta=\overline\Theta^k(V,x_0)\,,$$ 
then $\theta$ is called the $k$--dimensional density of $V$
at $x_0$ (or simply $k$--density at $x_o$).
\begin{remark} Recall \cite{federer} that if $V\subset\R^n$ is a
  semialgebraic subset, then the above two limits are equal and the
  $k$--density of $V$ is well defined for any point of $\R^n$.
\end{remark}
\begin{definition}
  Let $X\subset\R^n$ be a $k$-dimensional semialgebraic set and let
$x_0\in X$ be a point such that the $k$--density of $X$ at
$x_0$ is positive. A $(k-1)$--dimensional rectifiable subset
$Y\subset X$ such that $x_0\in Y$ is called a \emph{separating set of
  $X$ at $x_0$} if (see Fig.\ \ref{Figure 1})
\begin{itemize}
\item for some small $\epsilon>0$ the subset
  $\bigl(\epsilon B(x_0)\cap X\bigr)\setminus Y$ has at least two connected
  components $A$ and $B$,
\item the $(k-1)$--density of $Y$ at $x_0$ is  zero,
\item the inferior $k$--densities of $A$ and $B$ at $x_0$ are nonzero.
\end{itemize}
\end{definition}
\begin{proposition}[Lipschitz invariance of separating sets] Let $X$
  and $Z$ be two real semialgebraic sets. If there exists a
  bi-Lipschitz homeomorphism of germs $F\co (X,x_0)\rightarrow
  (Z,z_0)$ with respect to the inner metric, then $X$ has a separating
  set at $x_0\in X$ if and only if $Z$ has a separating set at $z_0\in
  Z$.
\end{proposition}
\begin{proof} The result would be immediate if separating sets were
  defined in terms of the inner metrics on $X$ and $Z$. So we must
  show that separating sets can be defined this way.

  Let $X\subset\R^n$ be a connected semialgebraic subset. Consider the
  set $X$ equipped with the inner metric and with the Hausdorff
  measure $\mathcal{H}_X^k$ associated to this metric. Let $V\subset
  X$ be a $k$-dimensional rectifiable subset. We define the inner
  inferior and superior densities of $V$ at $x_0\in X$ with respect to inner
  metric on $X$ as follows:
$$\underline\Theta^k(X,V,x_0)=\lim_{\epsilon\to 0^+}
\inf\frac{\mathcal{H}_X^k(V\cap \epsilon B_X(x_0))}{\eta\epsilon^k}$$
and
$$\overline\Theta^k(X,V,x_0)=\lim_{\epsilon\to
  0^+}\sup\frac{\mathcal{H}_X^k(V\cap
  \epsilon B_X(x_0))}{\eta\epsilon^k}\,,$$ 
where $\epsilon B_X(x_0)$
denotes the closed ball in $X$ (with respect to the inner metric) of
radius $\epsilon$ centered at $x_0$. 
The fact that separating sets can be defined using the inner metric
now follows from the following proposition, completing the
proof.\end{proof}
\begin{proposition}\label{proposition1} Let $X\subset\R^n$ be a
  semialgebraic connected subset. Let $V\subset X$ be a
  $k$-dimensional rectifiable subset and $x_0\in X$. Then, there exist
  two positive constants $\kappa_1$ and $\kappa_2$ such that:
$$\kappa_1\underline\Theta^k(X,V,x_0)\leq 
\underline\Theta^k(V,x_0)\leq \kappa_2\underline\Theta^k(X,V,x_0)$$
and
$$\kappa_1\overline{\Theta^k}(X,V,x_0)\leq 
\overline\Theta^k(V,x_0)\leq \kappa_2\overline\Theta^k(X,V,x_0)\,.$$
\end{proposition}
\begin{proof}
This follows immediately from the Kurdyka's ``Pancake
Theorem'' (\cite{Ku}, \cite{BM}) which says that if 
$X\subset\R^n$ is a semialgebraic subset then there exists a
  finite semialgebraic partition $X=\bigcup_{i=1}^{l} X_i$ such that
  each $X_i$ is a semialgebraic connected set whose inner metric
  and Euclidean metric are bi-Lipschitz equivalent.
\end{proof}

The following proposition shows, among other things, that the germ of
an isolated complex singularity which has a separating set cannot be
\emph{metrically conical}, i.e., bi-Lipschitz homeomorphic to the
metric cone on its link.
\begin{proposition}\label{cone} 
  Let $(X, x_0)$ be a ($n+1$)-dimensional metric cone whose base is a
  compact connected Lipschitz manifold (possibly with
  boundary). Then, $X$ does not have a separating set at $x_0$.
\end{proposition}
\begin{proof} Let $M$ be an $n$-dimensional compact connected
  Lipschitz manifold with boundary. For convenience of exposition we
  will suppose that $M$ is a subset of the Euclidean sphere
  $S^{k-1}\in\R^k$ centered at $0$ and with radius $1$ and $X$ the
  cone over $M$ with vertex at the origin $0\in\R^k$.  Suppose that
  $Y\subset X$ is a separating set, so $X\setminus Y=A\cup B$ with $A$
  and $B$ open in $X\setminus Y$; the $n$--density of $Y$ at $0$ is
  equal to zero and the inferior $(n+1)$--densities of $A$ and $B$ at
  $0$ are unequal to zero. In particular, there exists $\xi >0$ such
  that these inferior densities of $A$ and $B$ at $0$ are bigger than
  $\xi$. For each $t>0$, let $\rho_t\co X\cap tD^k\rightarrow X$ be
  the map $\rho_t(x)=\frac{1}{t}x$, where $tD^k$ is the ball about
  $0\in \R^k$ of radius $t$.  Denote $Y_t=\rho_t(Y\cap tD^k)$,
  $A_t=\rho_t(A\cap tD^k)$ and $B_t=\rho_t(B\cap tD^k)$. Since the
  $n$--density of $Y$ at $0$ is equal to zero, we have:
$$\lim_{t\to
    0^+}\mathcal{H}^{n}(Y_t)=0\,.$$ 
Also, since the inferior
  densities of $A$ and $B$ at $0$ are bigger than $\xi$, we have that
  $\mathcal{H}^{n+1}(A_t)>\xi$ and $\mathcal{H}^{n+1}(B_t)>\xi$ for all
  sufficiently small $t>0$. 

  Let $r$ be a radius such that $X\cap rD^k$ has volume $\le\xi/2$ and
  denote by $X'$, $A'_t$, $B'_t$, $Y'_t$ the result of removing from
  $X$, $A_t$, $B_t$, $Y_t$ the intersection with the interior of the
  ball $rB^k$. Then $X'$ is a Lipschitz $(n+1)$--manifold (with
  boundary), $A'_t$ and $B'_t$ subsets of $(n+1)$--measure $>\xi/2$
  separated by $Y_t$ of arbitrarily small $n$--measure.

The following lemma then gives the contradiction to
complete the proof.
\end{proof}
\begin{lemma}\label{lipschitz_manifold} Let $X'$ be a $(n+1)$-dimensional
  compact connected  Lipschitz manifold with boundary. Then, for
  any $\xi>0$ there exists $\epsilon>0$ such that if $Y'\subset X'$ is a
  $n$-dimensional rectifiable subset with
  $\mathcal{H}^{n}(Y')<\epsilon$, then $X'\setminus Y'$ has a connected
  component $A$ of $(n+1)$--measure exceeding
  $\mathcal{H}^{n+1}(X')-\xi/2$ (so any remaining components have total measure
  $<\xi/2$).
\end{lemma}
\begin{proof}
  If $X'$ is bi-Lipschitz homeomorphic to a ball then this follows
  from standard isoperimetric results: for a ball the isoperimetric
  problem is solved by spherical caps normal to the boundary (Burago
  and Maz'ja \cite{burago-mazja} p.\ 54, see also Hutchins
  \cite{hutchins}). Since the isoperimetric problem is generally
  formulated in terms of currents, one needs also that the mass of the
  current boundary of a region is less than or equal to the Hausdorff
  measure of the topological boundary (\cite{federer} 4.5.6 or
  \cite{morgan} Section 12.2).

  Let $\{ T_i \}_{i=1}^{m}$ be a cover of $X'$ by subsets which are
  bi-Lipschitz homeomorphic to balls and such that
$$T_i\cap T_j\neq\emptyset ~\Rightarrow ~\mathcal{H}^{n+1}(T_i\cap
T_j)>0.$$ Without loss of generality we may assume 
$$\xi/m<\min\{\mathcal{H}^{n+1}(T_i\cap T_j)~|~T_i\cap T_j\ne
\emptyset\}\,.$$

Since $T_i$ is bi-Lipschitz homeomorphic to a ball there
exists $\epsilon_i$ satisfying the conclusion of this lemma for
$\xi/m$.  Let $\epsilon=\min(\epsilon_1,\dots,\epsilon_m)$. So if $Y'\subset X'$ is an $n$--dimensional
rectifiable subset such that $\mathcal{H}^n(Y')<\epsilon$, then
for each $i$ the largest component $A_i$ of $T_i\setminus Y'$ has
complement $B_i$ of measure $<\xi/2m$.

We claim $\bigcup_{i=1}^mA_i$ is connected. It suffices to show that $$T_i\cap
T_j\ne \emptyset~\Rightarrow~A_i\cap A_j\ne\emptyset\,.$$ So suppose
$T_i\cap T_j\ne \emptyset$. Then $B_i\cup B_j$ has measure less than
$\xi$, which is less than $\mathcal H^n(T_i\cap T_j)$, so $T_i\cap
T_j\not\subset B_i\cup B_j$. This is equivalent to $A_i\cap A_j\ne \emptyset$.

Thus there exists a connected component $A$ of $X'\setminus Y'$ which
contains $\bigcup_{i=1}^mA_i$. Its complement $B$ is
a subset of 
$\bigcup_{i=1}^m B_i$ and thus has measure less than $\xi/2$.
\end{proof}

\section{Separating sets in normal surface singularities}

\begin{theorem}\label{th:sep} Let $X\subset\C^3$ be a weighted homogeneous
  algebraic surface with respect to the weights $w_1\geq w_2>w_3$ and
  with an isolated singularity at $0$. If
  $\bigl(X\setminus\{0\}\bigr)\cap\{z=0\}$ is not connected, then $X$
  has a separating set at $0$.
\end{theorem}
\begin{example}
  This theorem applies to the Brieskorn singularity $$X(p,q,r):=\{(x,y,z)\in
  \C^3~|~x^p+y^q+z^r=0\}$$ 
if $p\le q<r$ and $\operatorname{gcd}(p,q)>1$. In particular it is 
not metrically conical. This was known for a different reason by
\cite{BFN1}: a weighted homogeneous surface singularity (not
necessarily hypersurface) whose two lowest weights are distinct is not
metrically conical.
\end{example}
\begin{proof}[Proof of Theorem \ref{th:sep}]
  Take $\epsilon$ small enough that the intersection $V\cap\epsilon
  S^5$ is transverse and gives the singularity link.  Let $\tilde{A}$,
  $\tilde{B}\subset \bigl(V\cap \epsilon S^5\bigr)\cap\{ z=0\}$ be two
  semialgebraic closed subsets such that
  $\tilde{A}\cap\tilde{B}=\emptyset$. Let $\tilde{M}$ be the conflict
  set of $\tilde{A}$ and $\tilde{B}$ on $\epsilon S^5$, i.e.,
$$\tilde M=\{p\in\epsilon S^5 ~|~ d(p,\tilde A)=d(p,\tilde B)\}\,,$$
where $d(\cdot,\cdot)$ is the standard metric on $\epsilon S^5$
(euclidean metric in $\C^3$ gives the same set).  Clearly, $\tilde{M}$
is a compact semialgebraic subset and there exists $\delta>0$ such
that $d(\tilde{M},\{ z=0\})>\delta$. Let $M=\C^*\tilde M\cup\{0\}$
(the closure of the union of $\C^*$--orbits through $\tilde{M}$). Note
that the $\C^*$--action restricts to a unitary action of $S^1$, so the
construction of $\tilde M$ is invariant under the $S^1$--action, so
$M=\R^*\tilde M$, and is therefore $3$--dimensional. It is
semi-algebraic by the Tarski-Seidenberg theorem.  We will use the
weighted homogeneous property of $M$ to show $\dim(T_0M)\leq 2$, where
$T_0M$ denotes the tangent cone of $M$ at $0$, from which will follow
that $M$ has zero $3$--density. In fact, we
will show that $T_0M\subset \{x=0,y=0\}$.
  
Let $T\colon \tilde{M} \times [0,+\infty)\rightarrow
M$ be defined by:
$$T((x,y,z),t)=(t^{\frac{w_1}{w_3}}x,t^{\frac{w_2}{w_3}}y,tz).$$
Clearly, the restriction $T|_{ \tilde{M} \times (0,+\infty)}\colon
\tilde{M} \times (0,+\infty)\rightarrow M\setminus\{0\}$ is a
bijective semialgebraic map. Let $\gamma\colon [0,\epsilon)\rightarrow
M$ be a semianalytic arc; $\gamma(0)=0$ and $\gamma'(0)\neq 0$. We
consider $\phi(s)=T^{-1}(\gamma(s))$ for all $s\neq 0$.  Since $\phi$
is a semialgebraic map and $M$ is compact, 
$\displaystyle\lim_{s\to 0}\phi(s)$ exists and belongs to $M\times
\{0\}.$ For the same reason, $\displaystyle\lim_{s\to
  0}\phi^{\prime}(s)$ also exists and is nonzero. Therefore, the arc
$\phi$ can be extended to $\phi\colon [0,\epsilon)\rightarrow
\tilde{M} \times [0,+\infty)$ such that $\phi(0)\in \tilde{M}\times
\{0\}$ and $\phi'(0)$ exists and is nonzero. We can take the
$[0,\infty)$ component of $\phi$ as parameter and write
$\phi(t)=((x(t),y(t),z(t),t)$.  Then $\gamma(t)=(t^{w_1/w_3}x(t),
t^{w_2/w_3}y(t), tz(t))$, so
\begin{eqnarray*}
  \lim_{t\to 0^+}\frac{\gamma(t)}{t} &=& 
\left(\lim_{t\to 0}\frac{t^{\frac{w_1}{w_3}}}{t}x(t)\,,~
\lim_{t\to 0}\frac{t^{\frac{w_2}{w_3}}}{t}y(t)\,,~\lim_{t\to 0}z(t)\right) \\
  &=& (0,0,z(0))\,.
\end{eqnarray*}
This is a nonzero vector in the set $\{x=0,y=0\}$, so we obtain
that $$T_0\subset \{x=0,y=0\}.$$
  
Since $M$ is a $3$-dimensional semialgebraic set and $\dim(T_0M)\leq
2$, we obtain that the $3$-dimensional density of $M$ at $0$ is equal
to zero  (\cite{KR}).

Now, we have the following decomposition: $$V\setminus M=A\cup B\,,$$
where $\tilde{A}\subset A$, $\tilde{B}\subset B$, $A$ and $B$ are
$\C^*$--invariant and $A\cap B=\emptyset$. Since $A$ and $B$ are
semialgebraic sets, the $4$--densities $\operatorname{density}_4(A,0)$
and $\operatorname{density}_4(B,0)$ are defined. We will show that
these densities are nonzero. It is enough to prove that
$\dim_\R(T_0A)=4$ and $\dim_\R(T_0B)=4$. Let $\Gamma\subset A$ be a
connected component of $A\cap\{ z=0\}$. Note that
$\bar\Gamma=\Gamma\cup\{0\}$ is a complex algebraic curve.
We will show that $T_0A$ contains the set $\{(x,y,v)~|~ (x,y,0)\in
\bar\Gamma, v\in \C\}$ if $w_1=w_2$ (note that $\bar\Gamma$ is the line
through $(x,y,0)$ in this
case) or either the $y$--$z$ or the $x$--$z$ plane if $w_1<w_2$.

Given a smooth point $(x,y,0)\in \Gamma$ and $v\in \C$, we may choose
a smooth arc $\gamma\colon [0,\epsilon)\to A$ of the form
$\gamma(t)=(\gamma_1(t),\gamma_2(t), t^m\gamma_3(t))$ with
$(\gamma_1(0),\gamma_2(0))=(x,y)$ and $\gamma_3(0)=v$. Then, using the
$\R^*$--action, we transform this arc to the arc
$\phi(t)=t^j\gamma(t)$ with $j$ chosen so $jw_3+m=jw_2$. Now
$\phi(t)=(t^{jw_1}\gamma_1(t), t^{jw_2}\gamma_2(t),
t^{jw_2}\gamma_3(t))$ is a path in $A$ starting at the origin. Its
tangent vector $\rho$ at $t=0$,
$$\rho=\lim_{t\to 0+}\frac{\phi(t)}{t^{jw_2}}\,,$$
is $\rho=(x,y,v)$ if $w_1=w_2$ and $\rho=(0,y,v)$ if $w_1>w_2$. If
$w_1>w_2$ and $y=0$ then the same argument, but with $j$ chosen with
$jw_3+m=jw_1$, gives $\rho=(x,0,v)$. This proves our claim and
completes the proof that $T_0A$ has real dimension $4$. The proof for
$T_0B$ is the same.  
\end{proof}

\section{The Brian\c con-Speder example}

For each $t\in\C$, let $X_t=\{(x,y,z)\in\C^3 ~|~
x^5+z^{15}+y^7z+txy^6=0 \}$. This $X_t$ is weighted homogeneous
with respect to weights $(3,2,1)$ and has an isolated
singularity at $0\in\C^3$. 

\begin{theorem}\label{bs_1} $X_t$ has a separating set at $0$ if $t\ne
  0$ but does not have a separating set at $0$ if $t=0$.
\end{theorem}
\begin{proof}
  Note that for $t\ne 0$ Theorem \ref{th:sep} applies, so $X_t$ has a
  separating set.  So from now on we take $t=0$. Denote $X:=X_0$.  In
  the following, for each sufficiently small $\epsilon>0$, we use the
  notation $$X^{\epsilon}=\{ (x,y,z)\in X ~|~ \epsilon |y|\leq |z|
  \leq \frac{1}{\epsilon}|y| \}.$$
We need a lemma.
\begin{lemma}\label{lemma1} $X^{\epsilon}$ is
  metrically conical at the origin with connected link. 
\end{lemma}
\begin{proof}
Note that the lemma makes a statement about the germ of $X^\epsilon$
at the origin. We will restrict to the part of $X^\epsilon$ that
lies in a suitable closed neighborhood of the origin.

Let $P\co \C^3\rightarrow\C^2$ be the orthogonal projection
$P(x,y,z)=(y,z)$. The restriction $P_X$ of $P$ to $X$ is a $5$-fold
cyclic branched covering map branched along
$\{(y,z)~|~z^{15}+y^7z=0\}$. This is the union of the $y$--axis in
$\C^2$ and the seven curves $y=\zeta z^2$ for $\zeta$ a $7$--th root
of unity. These seven curves are tangent to the $z$--axis.

Let $$C^{\epsilon}=\{ (y,z)\in \C^2 ~|~ \epsilon |y|\leq |z| \leq
\frac{1}{\epsilon}|y| \} .$$ Notice that no part of the branch locus
of $P_X$ with $|z|<\epsilon$ is in $C^{\epsilon}$. In particular, if
$D$ is a disk in $\C^2$ of radius $<\epsilon$ around $0$, then the map
$P_X$ restricted to $X^\epsilon$ has no branching over this disk. We
choose the radius of $D$ to be $\epsilon/2$ and denote by $Y$ the part
of part of $X^\epsilon$ whose image lies inside this disk.
Then $Y$ is a covering of $C^\epsilon\cap D$, and to complete the
proof of the lemma we must show it is a connected covering space and
that the covering map is bi-Lipschitz.

Since it is a Galois covering with group $\Z/5$, to show it is a
connected cover it suffices to show that there is a closed curve in
$C^\epsilon\cap D$ which does not lift to a closed curve in
$Y$. Choose a small constant $c\le \epsilon/4$ and consider the curve
$\gamma\colon[0,1]\to C^\epsilon\cap D$ given by $\gamma(t)=(ce^{2\pi
  it}, c)$. A lift to $Y$ has $x$--coordinate $(c^{15}+c^8e^{14\pi i
  t})^{1/5}$, which starts close to $c^{8/5}$ (at $t=0$) and ends close
to $c^{8/5}e^{(14/5)\pi i}$ (at $t=1$), so it is not a closed curve.

To show that the covering map is bi-Lipschitz, we note that locally
$Y$ is the graph of the implicit function $(y,z)\mapsto x$ given by
the equation $x^5+z^{15}+y^7z=0$, so it suffices to show that the
derivatives of this implicit function are bounded. Implicit
differentiation gives
$$\frac{\partial x}{\partial  y}
=-\frac{7y^6z}{5x^4}\,,\quad\frac{\partial x}{\partial
  z}=-\frac{15z^{14}+y^7}{5x^4}\,.$$
 It is easy to see that there exists
$\lambda >0$ such that
$$|15z^{14}+ y^7|\leq \lambda|z|^4, \  
|y^7|\leq\lambda |z^{14}+y^{7}| \ \mbox{and} \
|15z^{14}+y^7|\leq\lambda |z^{14}+y^7|,$$ 
for all $(y,z)\in
C^{\epsilon}\cap D.$ We then get
$$\left|\frac{\partial x}{\partial
    y}\right|^5=\frac{7^5|y^{30}z^5|}{5^5|z^{14}+y^7|^4|z|^4}\le\frac{7^5}{5^5}\lambda^4|y^2z|<\frac{7^5\lambda^4\epsilon^3}{5^52^3}$$
and
$$\left|\frac{\partial x}{\partial
    z}\right|^5=\frac{|15z^{14}+y^7|^5}{5^5|z^{14}+y^7|^4|z|^4}\le \frac{\lambda^5}{5^5}\,,$$
completing the proof.
\end{proof}

We now complete the proof of Theorem \ref{bs_1}. Let us suppose that
$X$ has a separating set. Let $A,B,Y\subset X$ be subsets satisfying:
\begin{itemize}
\item for some small $\epsilon>0$ the subset $[\epsilon B(x_0)\cap
  X]\setminus Y$ is the union of relatively open subsets  $A$ and $B$,
\item the $3$-dimensional density of $Y$ at $0$
 is equal to zero,
\item the $4$-dimensional inferior densities of $A$ and $B$ at $0$
 are unequal to zero.
\end{itemize}
Set $$N^{\epsilon}=\{(x,y,z)\in\C^3 ~|~ |z|\leq\epsilon |y| \ \mbox{or}
\ |y|\leq\epsilon |z| \}.$$ For each subset $H\subset\C^3$ we
denote $$H^{\epsilon}=H\cap [\C^3\setminus N^{\epsilon}].$$

In this step, it is valuable to observe that there exists a positive constant
$K$ (independent of $\epsilon$) such that
\begin{equation}\label{volume_estimative}
\mathcal{H}^4(X\cap N^{\epsilon}\cap B(0,r)) \leq K  \epsilon r^4
\end{equation}
for all $0< r\leq 1$ (see, e.g., Comte-Yomdin \cite{CY}, chapter
5). By definition, the $4$-dimensional inferior density of $A$ at $0$
is equal to
$$\liminf_{r\to 0^+}\left(\frac{\mathcal{H}^4(A^{\epsilon}\cap
    B(0,r))}{r^4}
+\frac{\mathcal{H}^4(A\cap N^{\epsilon}\cap B(0,r))}{r^4}\right)$$
Then, if $\epsilon>0$ is sufficiently small, we can use inequality
\eqref{volume_estimative} in order to show that the $4$-dimensional
inferior density of $A^{\epsilon}$ is a positive number. In a similar
way, we can show that if $\epsilon>0$ is sufficiently small, then the
$4$-dimensional inferior density of $B^{\epsilon}$ at $0$ is a
positive number. These facts are enough to conclude that
$Y^{\epsilon}$ is a separating set of $X^{\epsilon}$.  But in view of 
 Lemma \ref{lemma1} this contradicts Proposition \ref{cone}.
\end{proof}

\end{document}